\documentclass [11pt] {article}

\usepackage {amssymb}

\title {A Remark on Formal KMS States in Deformation Quantization}

\author {{\bf
          Martin 
          Bordemann\thanks{Martin.Bordemann@physik.uni-freiburg.de}~,
          \addtocounter{footnote}{2}
          Hartmann 
          R{\"o}mer\thanks{Hartmann.Roemer@physik.uni-freiburg.de}~,
          Stefan 
          Waldmann\thanks{Stefan.Waldmann@physik.uni-freiburg.de}
         } \\[3mm]
         Fakult\"at f\"ur Physik\\Universit\"at Freiburg \\
         Hermann-Herder-Str. 3 \\
         79104 Freiburg i.~Br., F.~R.~G \\[3mm]
        }

\date{FR-THEP-98/1 \\[1mm]
      January 1998 \\[5mm]}

\newcommand {\im} {{\bf i}}
                                  
\newcommand {\BEQ} [1] {\begin {equation} \label {#1}}
\newcommand {\EEQ} {\end {equation}}

\newcommand {\Lie} {\mathcal L}
\newcommand {\cc} [1] {\overline {{#1}}}
\newcommand {\supp} {{\rm supp}}

\newcommand {\id} {{\sf id}}
\newcommand {\tr} {{\sf tr}}
\newcommand {\ad} {{\rm ad}}

\newcommand {\Exp} {{\rm Exp}}

\newenvironment {proof}{\small {\sc Proof:}}{{\hspace*{\fill} $\square$}}

\newtheorem {lemma} {Lemma} [section]
\newtheorem {proposition} [lemma] {Proposition}
\newtheorem {theorem} [lemma] {Theorem}
\newtheorem {corollary} [lemma] {Corollary}
\newtheorem {definition}[lemma] {Definition}

\begin {document}

\maketitle

\begin {abstract}
In the framework of deformation quantization we define formal KMS states
on the deformed algebra of power series of functions with compact 
support in phase space as $\mathbb C[[\lambda]]$-linear functionals
obeying a formal variant of the usual KMS condition
known in the theory of $C^*$-algebras. We show that for each temperature
KMS states always exist and are up to a normalization equal to the 
trace of the argument multiplied by a formal
analogue of the usual Boltzmann factor, a certain formal star exponential.
\end {abstract}

\newpage

\section {Introduction}
\label {IntroSec}

The concept of deformation quantization has been set up in \cite{BFFLS78}
and the existence of formal associative deformations of the pointwise
multiplication in the space of all complex-valued smooth functions on
a symplectic manifold, the so-called star products, has been proved in 
\cite{DL83}. The deformed algebra can be seen as a module over 
the formal power series ring $\mathbb C[[\lambda]]$ where the deformation 
parameter $\lambda$ corresponds to Planck's constant $\hbar$ and the 
constructions can be made such that the pointwise complex conjugation 
becomes an antilinear involutive antiautomorphism of the deformed algebra.
Using the natural ring ordering in the subring 
$\mathbb R[[\lambda]]$ of real power series
it is possible to define formal positive $\mathbb C[[\lambda]]$-linear
functionals on the deformed algebra and to imitate the GNS construction
known in the theory of complex $C^*$-algebras (see \cite{BW96b,BW97b}) 
which gives a notion of formal states in the theory of 
deformation quantization yielding physically reasonable representations
such as the Schr\"odinger picture or the WKB expansion for cotangent
bundles (see \cite{BNW97a,BNW97b}).

Having a notion of formal states it is natural to consider problems
of quantum statistical physics in this light. In the modern approach
based on the quantum observable algebra (taken to be a complex 
$C^*$-algebra) the analogue of a Gibbs state of inverse temperature 
$\beta$ is a positive linear functional $\mu$ on the algebra obeying
the so-called KMS condition (see for example the books by Bratteli-Robinson
\cite{BR81}, Haag \cite{Haa93}, or Connes \cite{Con94} or 
Section \ref{KMSSec} of this paper for a precise definition).
Originally, the KMS condition appeared as a boundary condition for
complex times for thermal Green functions 
in the papers of Kubo and Martin \& Schwinger
(see \cite{Kub57} and \cite[p.1357, p.1359]{MS59})
and was cast into the $C^*$-algebra language by Haag, Hugenholtz, and 
Winnink \cite {HHW67}. This condition
proved to be rather useful in the development of the statistical theory 
based on $C^*$-algebras, and it is believed that the nonuniqueness
of KMS states for a certain temperature is related to the existence of
several different thermodynamic phases
(see \cite{HHW67}, or \cite[p.213]{Haa93}, or \cite[p.41]{Con94} for a 
discussion). 

Beside the usual approaches in quantum field theory investigations 
in this directions in the setting of classical mechanics of infinitely 
many degrees of freedom has been made in e.~g. \cite {AGGL76,GP76,GV75} 
where the situation of infinitely many particles moving in flat 
$\mathbb R^n$ is considered by 
using sequences of coordinates and momenta for the particles and 
measure theoretical techniques to describe the KMS states.

More than ten years ago \cite{BFLS84,BL85} have already given a
treatment of the KMS condition in the framework of deformation 
quantization: the inverse temperature $\beta$ is incorporated in
the deformed algebra by an equivalence transformation (having zeroth
order term {\em not} equal to the identity) which is a
left multiplication by an invertible, $\beta$-dependent function
(such as an analogue of the Boltzmann factor). Rigidity and equivalence
of such $\beta$-dependent star products have been further discussed
in \cite{BL85}. The connection to the KMS condition is made by assuming
the existence of some complex topological subalgebra $\mathcal A$ of the
deformed algebra and the existence of a {\em complex} continuous linear 
functional $\mu$ on $\mathcal A$ such that the KMS condition for
a Hamiltonian $H$ makes sense, and by deriving a condition on the
$\beta$-dependent star product (see \cite{BFLS84}, Section 3, p.~490,
in particular eqn (3.3), and eqn (3.10) on p.~492).

More recently the classical KMS condition in the context of general
Poisson manifolds has been discussed in \cite {Wei97}: the starting point 
is a positive density $\mu$ on the manifold whose Lie derivative with 
respect to a Hamiltonian vector field gives rise to a unique vector field 
$\phi_\mu$ the so-called modular vector field which can be regarded as an 
infinitesimal version of the modular automorphisms in the Tomita-Takesaki 
theory of von Neumann algebras.

In this Letter we shall discuss the simplest case of finitely many degrees
of freedom: using the formal subalgebra $C_0^{\infty}(M)[[\lambda]]$ 
of series of smooth complex-valued functions having compact support 
in a connected symplectic manifold $M$ we show first
that for any Hamiltonian function $H$ the KMS-condition can be formulated
in terms of $\mathbb C[[\lambda]]$-linear functionals $\mu$ on
$C_0^{\infty}(M)[[\lambda]]$. Secondly, we prove ---without making any a priori
assumptions on the continuity of the functionals with respect to the standard 
locally convex topology--- that there is always a {\em unique} 
(up to normalizations in $\mathbb C[[\lambda]]$) formal KMS state 
$\mu$ on $C_0^{\infty}(M)[[\lambda]]$ given by
the following analogue of the Boltzmann factor
\BEQ {FirstEq}
    \mu (f) = c \, \tr \left(\Exp (-\beta H) * f \right),
\EEQ
where $\tr$ is a nonzero trace on $C_0^{\infty}(M)[[\lambda]]$ (a
$\mathbb C[[\lambda]]$-linear functional on $C_0^{\infty}(M)[[\lambda]]$ 
vanishing on commutators) and $\Exp (-\beta H)$ is the star-exponential
of $-\beta H$ where no formal convergence problem arises since there is
no $\frac{1}{\lambda}$ in front of $H$ in the exponent. 
In case the complex
conjugation is an antilinear antiautomorphism of the star product the
prefactor can be chosen such that the KMS states become formally positive.
Thirdly, we show that for $\beta \ne 0$ there is {\em no} nonzero KMS 
state in case the quantum time development is induced by a symplectic, 
but non-Hamiltonian vector field.

Assuming for a moment that phase transitions are related to the 
nonuniqueness of KMS states for a given inverse temperature $\beta$ 
we see that our result is physically reasonable insofar that phase 
transitions become mathematically visible only when some kind of 
thermodynamic limit is performed where particle 
number and configuration space volume are both sent to infinity while 
the average particle density is kept fixed: hence for 
finite-dimensional symplectic manifolds one would {\em not} expect 
phase transitions on physical grounds. For future investigations one 
would have to incorporate the
technically more involved formulation of thermodynamic limits (possibly
based on the work of \cite {AGGL76,GP76,GV75}) in the analysis and 
again look for formal KMS states in a more general infinite-dimensional 
setting.

The paper is organized as follows: Firstly we remember some basic facts 
concerning the notion of time development in deformation quantization 
as well as the notion of traces, i.~e. $\mathbb C[[\lambda]]$-linear 
functionals which vanish on star product commutators. Here we refer to the 
existence and uniqueness of traces established in \cite {NT95a} and give an 
alternative simple proof for the uniqueness of the traces which we shall 
need afterwards. In the following section we define formal KMS states 
after a short discussion of the KMS condition used in the context of 
$C^*$-algebras. Finally section \ref {ExistSec} contains the 
two main theorems: 
we prove the existence and uniqueness of KMS states for any star 
product on a connected symplectic manifold for any inverse temperatur 
$\beta$ with respect to the time development induced by an arbitrary 
Hamiltonian vector field. Moreover we show that for symplectic but 
non-Hamiltonian vector fields no KMS states for $\beta \ne 0$ exist.

\section {Some basic concepts: time evolution and traces}
\label {BasicSec}

In this section we shall briefly remember some basic facts concerning time 
evolution and traces as well as star exponentials in deformation 
quantization. Firstly we shall fix some notation: we consider a 
symplectic manifold 
$(M, \omega)$ and a symplectic vector field $X$. Then $i_X\omega = \alpha$ 
is a closed one-form and any closed one-form determines a symplectic 
vector field via this equation. By $\phi_t$ we denote the flow of $X$
where we assume for simplicity that $X$ has a complete flow. Then the 
classical time evolution of the observables, i.~e. the complex-valued
functions $C^\infty (M)$, with respect to $X$ is given by the pull-back 
$\phi_t^*: C^\infty (M) \to C^\infty (M)$ and for any initial 
condition $f \in C^\infty (M)$ the time evolution $f(t)$ through 
$f(0) = f$ is uniquely determined by 
\BEQ {HamiltonEqI}
    \frac{d}{dt} f(t) = \Lie_X f(t), \qquad f(0) = f
\EEQ
where $\Lie_X$ denotes the Lie derivative with respect to $X$. 
In the case where $\alpha = dH$ is exact the symplectic vector field $X$ 
is a Hamiltonian vector field with Hamiltonian function $H$ and 
(\ref {HamiltonEqI}) can be rewritten as
\BEQ {HamiltonEqII}
    \frac{d}{dt} f(t) = \left\{ f(t), H \right\}, \qquad f(0) = f 
\EEQ
where $\{ \cdot, \cdot \}$ denotes the Poisson bracket induced by 
$\omega$. Now in deformation quantization (see e.~g. \cite {BFFLS78})
the classical Poisson algebra $C^\infty (M)$ 
of observables is deformed into an associative star product algebra
$(C^\infty (M)[[\lambda]], *)$ where the star product $*$ is given by
the formal power series in the formal parameter $\lambda$
\BEQ {StarProd}
    f * g = \sum_{r=0}^\infty \lambda^r M_r (f, g)
\EEQ
with $M_0 (f, g) = fg$ and $M_1 (f, g) - M_1 (g, f) = \im \{f, g\}$
and all $M_r$ are bidifferential operators on $M$ vanishing on the constants 
for $r \ge 1$.
Here the deformation parameter $\lambda$ corresponds directly to 
Planck's constant $\hbar$ whence it is considered to be real, 
i.~e. we define $\cc \lambda := \lambda$. In the case of convergence 
$\lambda$ may be substituted by $\hbar \in \mathbb R$.
In the case of a Hamiltonian vector field $X$ the quantum analogue to 
(\ref {HamiltonEqII}) is given by Heisenberg's equation of motion
\BEQ {HeisenbergEqII}
    \frac{d}{dt} f(t) = \frac{\im}{\lambda} \ad(H) f(t) 
    \qquad 
    f(0) = f
\EEQ
where $\ad(H)g := H*g - g*H$ as usual and computing the first order 
in $\lambda$ of (\ref {HeisenbergEqII}) this can be viewed as 
deformation of (\ref {HamiltonEqII}). In the case where $X$ is only 
symplectic, i.~e. the corresponding one-form $\alpha$ is only closed but 
not exact one observes that Heisenberg's equation of motion can still be 
formulated: locally $\alpha = dH$ and using these locally defined 
Hamiltonians one observes that the locally defined map $\ad (H)$ 
only depends on $dH = \alpha$ and thus is indeed a globally defined 
map which we shall denote by $\delta_X$. Then $\delta_X$ is 
a derivation of the star product algebra. 
Fundamental for the following is the well-known existence of 
solutions $f(t) = A_t f$ of (\ref {HeisenbergEqII}) for any initial 
condition 
(see e.~g. \cite [Sec. 5.4] {Fed96}, \cite [App. B] {BNW97b}, 
\cite [Sec. 5]{BW96b}):
\begin {proposition} \label {TimeProp}
Let $(M, \omega)$ be a symplectic manifold and $X$ a symplectic vector 
field with complete flow $\phi_t$. Moreover let 
$*$ be a star product for $M$ then the Heisenberg equation of motion 
\BEQ {HeisenbergEqI}
    \frac{d}{dt} f(t) = \frac{\im}{\lambda} \delta_X f(t),
    \qquad 
    f(0) = f \in C^\infty (M)[[\lambda]] 
\EEQ
has a unique solution denoted by $f(t) = A_t f$ and the map 
$A_t: C^\infty (M)[[\lambda]] \to C^\infty (M)[[\lambda]]$ is a
$\mathbb C[[\lambda]]$-linear automorphism of $*$ 
and has the following properties:
\begin {enumerate}
\item $A_t = \phi^*_t \circ T_t$ where 
      $T_t = \id + \sum_{r=1}^\infty \lambda^r T^{(r)}_t$ and $T^{(r)}_t$ 
      is a differential operator vanishing on constants.
\item $A_t \circ \delta_X = \delta_X \circ A_t$ and $A_t$ is a 
      one-parameter group of automorphisms of the star product $*$.
\item If the complex conjugation is an antilinear anti-automorphism of 
      $*$, i.~e. $\cc{f * g} = \cc g * \cc f$ where 
      $\cc \lambda := \lambda$ then $A_t$ is a real automorphism, 
      i.~e. $\cc{A_t f} = A_t \cc f$.
\end {enumerate}
\end {proposition}      
In the following we shall often make use of a particular form of the star 
exponential \cite{BFFLS78} of a function $H \in C^\infty (M)$.
In our case the star exponential can be defined as a solution of a
differential equation which is shown to exist. In fact all relevant 
properties can be proved easily this way. We consider the differential 
equation
\BEQ {ExpDefEq}
    \frac{d}{d\beta} f(\beta) = H*f(\beta),
    \qquad 
    \beta \in \mathbb R .
\EEQ
\begin {lemma} \label {ExpLem}
Let $(M, \omega)$ be a sympectic manifold and $*$ a star product for $M$ 
and let $H \in C^\infty (M)$. Then there exists a unique solution 
$f(\beta)$ of (\ref {ExpDefEq}) in $C^\infty (M)[[\lambda]]$
with initial condition $f(0) = 1$. This 
solution is denoted by $\Exp(\beta H)$ and one has the following 
properties for all $\beta, \beta' \in \mathbb R$:
\begin {enumerate}
\item $\Exp (\beta H) = e^{\beta H}
      \left(1 + \sum_{r=1}^\infty \lambda^r g^{(r)}_\beta\right)$ where 
      $g^{(r)}_\beta \in C^\infty (M)$.
\item $\Exp (\beta H) * H = H * \Exp (\beta H)$ 
      and $\Exp (\beta H) * \Exp (\beta' H) = \Exp ((\beta + \beta')H)$.
\end {enumerate}
\end {lemma}
\begin {proof}
This lemma is proved by first factorizing 
$f(\beta) = e^{\beta H} g(\beta)$ and then rewriting the induced 
differential equation for $g(\beta)$ as integral equation which can be 
uniquely solved by recursion since the integral operator raises the 
degree in $\lambda$. Then the other properties 
easily follow using the uniqueness and (\ref {ExpDefEq}).
\end {proof}

Now we consider again a symplectic vector field $X$ and the corresponding 
derivation $\delta_X$. Since clearly 
$\delta_X = -\im\lambda \Lie_X + \cdots$ the map $\delta_X$ 
raises the $\lambda$-degree at least by one which implies that the series
\BEQ {ExpadalphaDef}
    e^{\beta\delta_X} 
    := \sum_{r=0}^\infty \frac{1}{r!} \left(\beta \delta_X\right)^r
\EEQ
is a well-defined formal power series of maps for $\beta \in \mathbb R$
and one easily shows that 
$e^{\beta \delta_X}$ is a one-parameter group of automorphisms 
of the star product. Moreover one has the following lemma:
\begin {lemma} \label {InnerLem}
Let $(M, \omega)$ be a symplectic manifold and $*$ a star product for $M$ 
and let $X$ be a symplectic vector field. Then the one-parameter group 
$e^{\beta\delta_X}$ of automorphisms of $*$ where 
$\beta \in \mathbb R$ is inner iff $i_X\omega = dH$ is exact 
and in this case for all $f \in C^\infty (M)[[\lambda]]$
\BEQ {ExpAdInner}
    e^{\beta\delta_X} (f) = \Exp (\beta H) * f * \Exp (-\beta H) .
\EEQ    
\end {lemma}
\begin {proof}
Let us first assume that $e^{\beta\delta_X}$ is an inner automorphism 
for some $\beta \ne 0$, i.~e. we assume that there exist elements 
$b = \sum_{r=0}^\infty \lambda^r b_r$ and 
$c = \sum_{r=0}^\infty \lambda^r c_r$
where $b_r, c_r \in C^\infty (M)$ (depending on $\beta$) such that 
\[
    e^{\beta \delta_X} (f) = b * f * c 
    \quad
    \mbox { and }
    \quad
    b*c = 1 = c*b .
\]
By straight forward computation of the first order in $\lambda$ of the 
relation $e^{\beta\delta_X} (f) - e^{-\beta\delta_X} (f) 
= b*f*c - c*f*b$ one obtains $\beta \Lie_X f = c_0 \{f, b_0\}$. Since 
$b*c=1$ one has $b_0 c_0 = 1$ and thus $b_0$ is a non-vanishing 
function on $M$. Now define
$H := \frac{1}{2\beta} \ln (b_0 \cc b_0)$ which is clearly a smooth 
function on $M$. We obtain $\{f, H\} = \Lie_X f$ which shows that $X$ 
is in fact Hamiltonian and thus $e^{\beta\delta_X}$ has only a 
chance to be inner if $i_X\omega = dH$ is exact. 
If on the other hand $i_X\omega = dH$ then (\ref {ExpAdInner}) is a 
simple computation using lemma \ref {ExpLem} and (\ref {ExpDefEq}).
\end {proof}

A last important structure needed in the following is the notion of a 
trace. Though in deformation quantization 
traces are usually considered in the setting of
formal Laurent series (e.~g. in \cite {Fed96,NT95a,CFS92}) which 
allows a more suitable normalization motivated by either physical 
reasons (`to get dimensions right') or by analogy to traces of 
pseudo-differential operators we shall stay for simplicity in the 
category of formal power series.
\begin {definition}
Let $(M, \omega)$ be a symplectic manifold and $*$ a star product for $M$. 
A $\mathbb C[[\lambda]]$-linear functional 
$\tr : C^\infty_0 (M)[[\lambda]] \to \mathbb C[[\lambda]]$ is called a 
trace iff $\tr (f*g - g*f) = 0$ for all 
$f, g \in C^\infty_0 (M)[[\lambda]]$.
\end {definition}
\begin {proposition} [Existence and uniqueness 
                      of traces \cite {NT95a,Fed96}]
\label {TraceProp}
Let $(M,\omega)$ be a connected and symplectic manifold and $*$ a 
star product for $M$. 
Then the set of traces forms a $\mathbb C[[\lambda]]$-module which is 
one-dimensional over $\mathbb C[[\lambda]]$.
\end {proposition}
For a proof of the existence we refer to e.~g. \cite {NT95a} where also 
the uniqueness up to normalization by elements in $\mathbb C[[\lambda]]$
is shown. For the existence of strongly closed star products as defined 
in \cite {CFS92} see \cite {OMY92}.
For later use we shall give here an elementary proof of the 
uniqueness since we need in particular the lowest (non-trivial) order of 
the traces. Though the following lemma should be well-known
we shall give for completeness a short proof since the result is crucial 
for the following:
\begin {lemma}
\begin{enumerate}
\item Let $\varphi \in C^\infty_0 (\mathbb R^n)$ ($n\geq 1$) be a smooth 
      complex valued 
      function with support contained in an open ball $B_R (0)$ 
      around $0$ with radius $R > 0$ such that 
      $\int_{\mathbb R^n} \varphi (x) d^n x = 0$
      then there exist functions $h_i \in C^\infty_0 (\mathbb R^n)$ with 
      $\supp h_i \subset B_R (0)$ for $i = 1, \ldots, n$ such that
      $\varphi = \sum_{i=1}^n \frac{\partial h_i}{\partial x^i}$.
\item Let $B_R (0) \subseteq \mathbb R^n$ be an open ball around $0$ with 
      radius $R > 0$ (where $R = \infty$ is also allowed) and let 
      $\mu: C^\infty_0 (B_R (0)) \to \mathbb C$ be a 
      $\mathbb C$-linear (not necessarily continuous) functional such 
      that $\mu \left( \frac{\partial f}{\partial x^i} \right) = 0$
      for all $i = 1, \ldots, n$ and $f \in C^\infty_0 (B_R (0))$ 
      then $\mu$ is a distribution and given by
      \[
          \mu (f) = c \int_{B_R (0)} f(x) \; d^n x
      \]
      with some constant $c \in \mathbb C$.    
\end{enumerate}
\end {lemma}
\begin {proof} 
For the first part the case $n=1$ is readily checked by noting that there
is a primitive $h_1$ of $\varphi$ having compact support. Assume that
$n\geq 2$. Choose three pairwise distinct concentric closed balls 
$B_1\subset B_2\subset B_3$ in $B_R(0)$ such that
${\rm supp}\varphi\subset B_1$. Embed $B_R(0)$ as an open subset in the
sphere $S^n$ and extend $\varphi$ to a smooth complex-valued function
on $S^n$ vanishing outside the embedded $B_R(0)$. We can assume that
the embedding is volume preserving. Hence the integral of the extended
$\varphi$ over $S^n$ (with some suitable volume $\mu$) is zero. Since
the $n$th de Rham cohomology group of $S^n$ is well-known to be 
one-dimensional it follows by the de Rham Theorem that $\varphi\mu$ 
is exact, hence equal to $d\alpha$ where $\alpha$ is some $n-1$-form 
on $S^n$. Now $\varphi$ vanishes
on the complement of the embedded ball $B_1$ in $S^n$ which is 
diffeomorphic to $\mathbb R^n$ hence $\alpha$ is closed on that set.
By the Poincar\'e Lemma there is an $n-2$-form $\beta$ on that subset
such that $\alpha=d\beta$ on that subset. Choose a nonnegative smooth
function $\chi$ on the sphere being zero on the embedded $B_2$ and $1$ 
on the complement of the interior of the embedded $B_3$ it follows
that $\alpha':=\alpha-d(\chi\beta)$ is a globally defined $n-2$-form
on the sphere with support in the embedded $B_3$ such that 
$\varphi\mu = d\alpha'$. Pulling this back to the ball $B_R(0)$ 
we get the desired functions $h_i$ by the components of the pulled-back 
$\alpha'$.
For the second part notice that part one shows that the 
linear space of smooth functions of compact support
in the ball generated by derivatives of such functions is of codimension
one hence all the linear functionals having this subspace in their kernel
must be multiples of the integral.
\end {proof}

Using the preceding Lemma and a standard partition-of-unity argument
we also have the
\begin {corollary} \label {PoissonNullCor}
Let $(M, \omega)$ be a connected symplectic manifold and let 
$\mu: C^\infty_0 (M) \to \mathbb C$ be a linear functional vanishing on 
Poisson brackets, i.~e. $\mu \left( \{f, g\}\right) = 0$ for all 
$f, g \in C^\infty_0 (M)$ then there exists a complex number 
$c \in \mathbb C$ such that 
\[
    \mu (f) = c \int_M f \; \Omega
\]
where $\Omega = \omega \wedge \cdots \wedge \omega$ is the symplectic 
volume form.       
\end {corollary}
Now we consider a non-trivial trace $\tr$ for 
a connected symplectic manifold $(M, \omega)$ with star product $*$. 
Firstly we remember that any $\mathbb C[[\lambda]]$-linear 
functional of $C^\infty_0 (M)[[\lambda]]$ 
can be written as 
\[
    \tr = \sum_{r=0}^\infty \lambda^r \mu_r
\]
with $\mu_r: C^\infty_0 (M) \to \mathbb C$ due to \cite [Prop. 2.1] {DL88}
and we can assume without restriction that $\mu_0 \ne 0$. 
Then the trace property of $\tr$ obviously implies that $\mu_0$ vanishes 
on Poisson brackets and hence there exists a complex number $c_0 \ne 0$ 
such that 
\BEQ {TrLowestOrder}
    \mu_0 (f) = c_0 \int_M f \; \Omega .
\EEQ
Now if $\tilde \tr$ is another trace for $*$ then 
$\tilde \mu_0 (f) = \tilde c_0 \int_M f \Omega$ and thus 
$\tilde \tr - \frac{\tilde c_0}{c_0} \tr$ is again a trace starting at least 
with order $\lambda^1$. Thus one can construct recursively 
a formal power series $c = \frac{\tilde c_0}{c_0} + \cdots$ 
such that $\tilde \tr = c\, \tr$ which proves the uniqueness of traces in 
the connected case up to normalization.

\section {The formal KMS condition in deformation quantization}
\label {KMSSec}

After these preliminaries we can now discuss the meaning of the KMS 
condition in deformation quantization which was first discussed in this 
context in \cite {BFLS84}. 
In our approach we try to stay completely in the formal category and 
avoid any assumptions about convergence of the formal power series.
Moreover we restrict ourselves to finite-dimensional phase spaces.

Firstly we shall shortly remember the well-known definition of KMS states 
used e.~g. in algebraic quantum field theory in the context of 
$C^*$-algebras (see e.~g. \cite {BR81,Haa93,Con94}). 
Here the observable algebra $\mathcal A$ is a net of local 
$C^*$-algebras with inductive limit topology and the time develompent 
operator $\alpha_t : \mathcal A \to \mathcal A$ is a one-parameter group 
of $^*$-automorphisms of $\mathcal A$ and a state $\mu$ 
(i.~e. a positive linear functional) of $\mathcal A$ is called a 
KMS state for the inverse 
temperature $\beta = \frac{1}{kT}$ (where $k$ is Boltzmann's constant and 
$T$ the absolute temperature) if for any two observables 
$a, b \in \mathcal A$ there exists a continous function 
$F_{ab} : S_\beta \to \mathbb C$  
which is holomorphic inside the strip 
$S_\beta := \{z \in \mathbb C \; | \; 0 \le \mbox{ Im } z \le \hbar\beta \}$
such that for any real $t$
\BEQ {KMSinQF}
    F_{ab} (t) = \mu (\alpha_t (a) b) 
    \quad
    \mbox { and }
    \quad
    F_{ab} (t + \im\hbar\beta) = \mu (b \alpha_t (a)) .
\EEQ
This formulation replaces in a mathematically reasonable way the more 
intuitive requirement that for any two observables $a, b \in \mathcal A$ 
the state $\mu$ should satisfy
\[
    \mu (\alpha_t (a) b) = \mu (b \alpha_{t + \im\hbar\beta} (a))
\]
which is obviously not well-defined in general since there is a priori no 
sense of the complexification of the time development operator $\alpha_t$ 
to define $\alpha_{t + \im\hbar\beta}$. Nevertheless we shall see that in 
deformation quantization there is indeed a reasonable notion of 
such a `complexification' which avoids the usage of the holomorphic 
functions $F_{ab}$ which is not suitable in the formal setting since we 
want to treat $\hbar$ as formal! The key ingredient is 
the following simple lemma which follows directly from 
proposition \ref {TimeProp} 
and the definition of $e^{\beta\delta_X}$:
\begin {lemma}
Let $(M, \omega)$ be a symplectic manifold with star product $*$ and let 
$X$ be a symplectic vector field on $M$ with 
complete flow and let $A_t$ be the 
corresponding time development operator. Then the map
\BEQ {ComplexAt}
    A_{t+ \im \lambda \beta} := A_t \circ e^{\beta \delta_X}
\EEQ
where $t, \beta \in \mathbb R$ is an automorphism of the star product 
$*$ and 
$A_{t+\im\lambda\beta} \circ A_{t' + \im\lambda\beta'} =
A_{t+t' + \im\lambda(\beta + \beta')}$ for all 
$t, t', \beta, \beta' \in \mathbb R$.
\end {lemma}
This seems to be a reasonable definition for the 
`complexification' of $A_t$ in this particular situation. Note that this 
would make no longer sense in general if we tried to define 
$A_{t+\im\beta}$ for $\beta \in \mathbb R$.

Now we can define formal KMS states in deformation quantization in 
the following way: firstly we remember that even in the formal 
setting there is a both 
mathematically and physically reasonable notion of positive linear 
functionals in the 
case where the star product satisfies $\cc{f*g} = \cc g * \cc f$ using the 
natural ring ordering of $\mathbb R[[\lambda]]$. Such positive linear 
functionals give raise to a formal GNS construction as defined in details 
in \cite {BW96b}. Hence it appears natural to consider only such star 
products and search for formal KMS states within these positive linear 
functionals. But it will turn out that the formal KMS condition will 
in fact essentially imply (in the connected case after suitable 
normalization) the positivity and hence we shall not proceed this way 
but state the following definition:
\begin {definition} \label {KMSDef}
Let $(M, \omega)$ be a symplectic manifold with star product $*$ and let 
$X$ be a symplectic vector field on $M$ and let 
$\mu : C^\infty_0 (M) [[\lambda]] \to \mathbb C[[\lambda]]$ 
be a $\mathbb C[[\lambda]]$-linear functional.
\begin {enumerate}
\item $\mu$ satisfies the {\bf static formal KMS condition} for the inverse 
      temperature $\beta \in \mathbb R$ with respect to $X$ 
      iff for all $f, g \in C^\infty_0 (M)[[\lambda]]$
      \BEQ {StaticKMS}
          \mu (f * g) = \mu (g * e^{\beta\delta_X} (f)) .
      \EEQ
\item If $X$ has complete flow then 
      $\mu$ satisfies the {\bf dynamic formal KMS condition} for the 
      inverse temperature
      $\beta \in \mathbb R$ with respect to the time development operator
      $A_t$ iff for all $f, g \in C^\infty_0 (M)[[\lambda]]$ 
      \BEQ {DynamicKMS}
          \mu \left(A_t(f) * g\right) 
          = \mu \left(g * A_{t + \im\lambda\beta} (f)\right) .
      \EEQ
\end {enumerate}
\end {definition}
Clearly $\mu =0$ satisfies the formal KMS condition trivially and if we 
consider $\mu \ne 0$ we can assume that the first non-trivial order of 
$\mu$ is the zeroth order in $\lambda$.
Evaluating the first non-trivial order in $\lambda$ 
of the KMS conditions (\ref {StaticKMS}) resp. (\ref {DynamicKMS}) 
using (\ref {ExpadalphaDef}) and 
proposition \ref {TimeProp} one obtains the well-known classical 
KMS conditions namely the {\em static classical KMS condition}
\BEQ {SCKMS}
    \mu_0 \left( \{f, g\} - \beta g \Lie_X f \right) = 0 
    \qquad
    \forall f, g \in C^\infty_0 (M)
\EEQ
and in the case where $X$ has complete flow the 
{\em dynamical classical KMS condition}
\BEQ {DCKMS}
    \mu_0 \left( \left\{ \phi^*_t f, g \right\} 
                 - \beta g \Lie_X \phi_t^* f \right) = 0
    \qquad
    \forall f, g \in C^\infty_0 (M),
    \;
    \forall t \in \mathbb R 
\EEQ
which where discussed in literature earlier in various ways 
(See e.~g. \cite {AGGL76,BFLS84,GP76,GV75,Wei97} and references therein).

In the case of a complete flow the dynamical KMS conditions (both quantum 
and classical) imply clearly the static ones by setting $t=0$. But since 
$A_t$ commutes with $\delta_X$ resp. $\phi^*_t$ commutes with 
$\Lie_X$ the static KMS conditions imply the dynamical ones by replacing 
$f$ by $A_t(f)$ resp. $\phi^*_t f$. Hence we shall only consider the 
static KMS condition in the following and drop the somehow technical 
assumption of complete flow.

\section {Existence and uniqueness of formal KMS states}
\label {ExistSec}

In the case of a Hamiltonian time development, i.~e. if $i_X\omega = dH$ with 
some Hamiltonian function $H \in C^\infty (M)$ the structure of the 
formal KMS states in the sense of definition \ref {KMSDef} is completely 
clarified by the following theorem:
\begin {theorem} \label {HamiltonianKMSTheo}
Let $(M, \omega)$ be a symplectic manifold with star product $*$ and let 
$H \in C^\infty (M)$ be a Hamiltonian function with corresponding 
Hamiltonian vector field $X$ and let $\beta \in \mathbb R$.
\begin {enumerate}
\item Let $\mu : C^\infty_0 (M)[[\lambda]] \to \mathbb C[[\lambda]]$ be a 
      $\mathbb C[[\lambda]]$-linear functional. Then $\mu$ satisfies the 
      static formal KMS condition (\ref {StaticKMS}) iff the functional
      $\tilde \mu (f) := \mu \left( \Exp (\beta H) * f \right)$
      is a trace for $*$.
\item If $M$ is connected then the set of static formal KMS states is 
      one-dimensional over $\mathbb C[[\lambda]]$ and any static formal
      KMS state $\mu$ can be obtained by
      \BEQ {KMSstates}
          \mu (f) = c \, \tr \left(\Exp (-\beta H) * f \right)
      \EEQ
      where $\tr$ is a non-trivial fixed choice of a trace for $*$ 
      starting with lowest order zero and $c \in \mathbb C[[\lambda]]$.
\item Let $\mu_0 : C^\infty_0 (M) \to \mathbb C$ be a 
      $\mathbb C$-linear functional then $\mu_0$ satisfies the static
      classical KMS condition iff the functional 
      $\tilde\mu_0 (f) := \mu_0 \left (e^{\beta H} f\right)$
      vanishes on Poisson brackets. Hence if $M$ is connected $\mu_0$ is 
      of the form
      \BEQ {CKMSState}
          \mu_0 (f) = c_0 \int_M e^{-\beta H} f \; \Omega
      \EEQ
      with some constant $c_0 \in \mathbb C$.
\end {enumerate}
\end {theorem}      
\begin {proof}
Part one is a simple and straight forward computation using 
lemma \ref {ExpLem} and \ref {InnerLem}. Then the second part follows 
immediately by proposition \ref {TraceProp}. The third part is shown the 
same way by computation and corollary \ref {PoissonNullCor}.
\end {proof}
      
Note that no continuity properties of $\mu_0$ had to be assumed for the 
classical part of the theorem. In fact the algebraic condition 
(\ref {SCKMS}) implies continuity of $\mu_0$ with respect to the 
standard locally convex topology of $C^\infty_0 (M)$ since clearly 
(\ref {CKMSState}) defines a continuous functional.

In the case when the time develoment is given by a symplectic but not 
Hamiltonian vector field no non-trivial formal KMS states exist:
\begin {theorem} \label {SymplecticKMSTheo}
Let $(M, \omega)$ be a connected symplectic manifold with star product 
$*$ and let $X$ be a symplectic vector field on $M$ and let 
$0 \ne \beta \in \mathbb R$. 
If $\mu$ is a static formal KMS state with respect to $X$ and 
inverse temperature $\beta$ then either $\mu = 0$ or 
$\alpha := i_X\omega = dH$ is exact.
\end {theorem}
\begin {proof}
Since the static formal KMS condition (\ref {StaticKMS}) implies 
the classical one we only have to show that the classical static KMS 
condition (\ref {SCKMS}) implies either $\mu_0 = 0$ or $\alpha = dH$. 
Now let $\mu_0 : C^\infty_0 (M) \to \mathbb C$ be a linear 
functional satisfying (\ref {SCKMS}) then we take an atlas on $M$ of 
contractable charts $\{U_i\}_{i\in I}$ and local functions 
$H_i \in C^\infty (U_i)$ such that $\alpha|U_i = dH_i$ for all $i \in I$.
Consider $U_{ij} := U_i \cap U_j \ne \emptyset$ and let 
$C_{ij} \in \mathbb R$ be the constants such that 
$H_i|U_{ij} = H_j|U_{ij} + C_{ij}$. Now define 
$\mu^{(i)}: C^\infty_0 (U_i) \to \mathbb C$ by $\mu^{(i)} (f) := \mu(f)$ 
then $\tilde\mu^{(i)} (f) := \mu(e^{\beta H_i} f)$ is well-defined for 
$f \in C^\infty_0 (U_i)$ for any $i \in I$ and vanishes on Poisson 
brackets. Hence there exist constants $C_i \in \mathbb C$ such that for 
any $f \in C^\infty_0 (U_i)$
\[
    \mu (f) = C_i \int_{U_i} e^{-\beta H_i} f \; \Omega
\]
due to theorem \ref {HamiltonianKMSTheo}.
Thus for $U_{ij} \ne \emptyset$ this implies by a standard continuity 
argument $C_i e^{-\beta H_i}|U_{ij} = C_j e^{-\beta H_j}|U_{ij}$. Now if 
$C_i = 0$ then for any other $j \in I$ we obtain $C_j =0$ since $M$ is 
connected and hence $\mu = 0$. If on the other hand $C_i \ne 0$ then
$C_j \ne 0$ and $\frac{C_i}{C_j} > 0$. Thus we obtain 
$H_i|U_{ij} = H_j|U_{ij} + \frac{1}{\beta} \ln \frac{C_i}{C_j}$ and thus
$C_{ij} = \frac{1}{\beta} \ln \frac{C_i}{C_j}$. Hence the constants 
$C_{ij}$ clearly satisfy the cocycle identity which implies that $\alpha$ 
is exact.
\end {proof}

Finally we shall consider the case where the star products satisfies 
$\cc{f*g} = \cc g * \cc f$ for all $f, g \in C^\infty (M)[[\lambda]]$ 
where we set as usual $\cc \lambda := \lambda$. Now let $\tr'$ be a trace 
for $*$ then this property of $*$ guarantees that 
$\tr := \tr' + \cc \tr'$ (where $\cc \tr (f) := \cc{\tr (\cc f)}$) is 
also a trace of $*$ with the additional 
property that this trace is {\em real} in the following sense:
\BEQ {RealTrace} 
    \tr (\cc f) = \cc {\tr (f)}.
\EEQ
In the connected case a real trace $\tr$ is either a positive 
linear functional, i.~e. $\tr (\cc f * f) \ge 0$ in the sense of the 
ring ordering of $\mathbb R[[\lambda]]$ 
(where $a = \sum_{r=k}^\infty \lambda^r a_r \in \mathbb R[[\lambda]]$ 
is called positive iff $a_k > 0$), or $-\tr$ is a positive linear 
functional:
\begin {lemma}
Let $(M, \omega)$ be a symplectic connected manifold and let $*$ be a 
star product for $M$ such that $\cc{f * g} = \cc g * \cc f$ and 
let $\tr$ be a real non-vanishing trace. Then either $\tr$ or $-\tr$ 
is a positive linear functional and the Gel'fand ideal 
$\mathcal J_\tr := 
\{ f \in C^\infty_0 (M)[[\lambda]] \; | \; \tr(\cc f * f) = 0\}$ 
is $\{0\}$.
\end {lemma}
\begin {proof}
Since the lowest order of $\tr$ is proportional to the integration over 
$M$ with volume form $\Omega$ and since $\tr$ is real 
\cite [Lemma 2] {BW96b} implies the lemma.
\end {proof}

This lemma implies that in the case where $\cc{f*g} = \cc g * \cc f$ 
a formal KMS state can by rescaled to obtain a positive formal KMS state. 
Hence the algebraic formal KMS condition (static or dynamic) 
implies positivity in deformation quantization.

\section* {Acknowledgements}

We would like to thank Moshe Flato and Daniel Sternheimer for 
discussions and for pointing out reference \cite {Wei97}.

\end {document}